\documentclass[10pt, reqno, a4wide]{article}
\usepackage{amssymb, amsmath, amsthm}
\usepackage{ifthen}
\usepackage[all, cmtip]{xy}
\usepackage{pdfpages}
\usepackage{graphics}
\usepackage{mathrsfs}

\newtheorem{leer}{}[section]
\newtheorem{thm}[leer]{Theorem}

\newtheorem{rema}[leer]{Remark} 
\newtheorem{prop}[leer]{Proposition}
\newtheorem{lemm}[leer]{Lemma}
\newtheorem{coro}[leer]{Corollary}
\newtheorem{defi}[leer]{Definition}

\newtheorem{intthm}{Theorem}

\newtheorem{intcoro}[intthm]{Corollary}

\newcommand{\ol}[1]{\overline{#1}}
\newcommand{\ul}[1]{\underline{#1}}

\newcommand{\plim}{\mathop{\varprojlim}\limits}
\newcommand{\ilim}{\mathop{\varinjlim}\limits}

\newcommand{\AAA}{\mathscr{A}}
\newcommand{\BBB}{\mathscr{B}}

\newcommand{\KKK}{\mathscr{K}}

\newcommand{\UUU}{\mathscr{U}}

\newcommand{\Aa}{{\mathbb{A}}}
\newcommand{\Bb}{{\mathbb{B}}}

\newcommand{\Ll}{{\mathbb{L}}}

\newcommand{\Nn}{{\mathbb{N}}}

\newcommand{\Qq}{{\mathbb{Q}}}

\newcommand{\Tt}{{\mathbb{T}}}

\newcommand{\Zz}{{\mathbb{Z}}}

\newcommand{\uG}{{\ul{G}}}

\newcommand{\uT}{{\ul{\Tt}}}

\newcommand{\Aut}{\mathrm{Aut}}

\newcommand{\GL}{\mathrm{GL}}

\newcommand{\Hom}{\mathrm{Hom}}
\newcommand{\End}{\mathrm{End}}

\newcommand{\Gal}{\mathrm{Gal}}

\newcommand{\coker}{\mathrm{coker}}
\newcommand{\im}{\mathrm{im}}
\newcommand{\chara}{\mathrm{char}}

\date{}

\newcommand{\Zar}{\mathrm{Zar}}
\newcommand{\Sm}{\mathrm{Sm}}
\renewcommand{\Bb}{Z}

\title{Local to global principles\\ for homomorphisms of abelian schemes}

\begin{document}
\author{by Wojciech Gajda and Sebastian Petersen\\
\\
{\small\sl dedicated  to Moshe Jarden with admiration}\\
{\small\em on the occasion of his 80th birthday}}

\maketitle

\abstract{}

\parindent 0cm

Let $A$ and $B$ be abelian varieties defined over  the function field $k(S)$ of a smooth algebraic variety $S/k.$ We establish criteria,  
in terms of restriction maps to subvarieties of $S,$ for existence of various important classes of $k(S)$-homomorphisms from $A$ to $B,$ e.g., 
for existence of $k(S)$-isogenies. Our main tools consist of Hilbertianity methods, Tate conjecture as proven by Tate, Zarhin and Faltings, and of the 
{\sl minuscule weights conjecture} of Zarhin in the case, when the base field is finite.   
 \endabstract



\section{Introduction}
Let $S$ be a smooth variety over a finitely generated 
field $k$ of arbitrary characteristic. Let $\AAA$ and $\BBB$ be $S$-abelian schemes with
generic fibers $A$ and $B$ (respectively) defined over the function field $k(S)$.\footnote{{\bf keywords}: Abelian variety, Abelian scheme, Galois representation\\ 
{\bf mathclass}: Primary 11G10; Secondary 14F20.} In this paper we consider existence of certain classes of  $k(S)$-homomorphisms from 
$A$ to $B$, e.g., $k(S)$-isogenies,  
and provide local criteria in terms of restriction maps to subvarieties of $S.$ 
Furthermore we study existence of abelian subvarieties of $A$ in a similar way. 
Our first main result is the following local to global principle.     
\medskip

\begin{intthm}[Thm. \ref{hilbertmain0}] \label{thmA}Let $S$ be a smooth variety over a finitely generated 
field $k$. Let $\AAA, \BBB$ be abelian schemes over $S$ with
generic fibers $A$ and $B,$ respectively. Let $U$ be a dense open subscheme of $S$. 
Let $m\in \{0,1,\cdots, \dim(S)\}$.  Assume that $k$ is infinite or that $m\ge 1$. Let $\kappa \in\Nn$. 

\begin{enumerate}
\item[(a)] The following are equivalent:
\begin{enumerate}
\item[(i)] There exists a $k(S)$-isogeny (resp. surjective homomorphism, resp. non-zero homomorphism, resp. homomorphism with $\kappa$-dimensional kernel)
$A\to B$.
\item[(ii)] For every $m$-dimensional smooth connected subscheme $T$ of $U$ there exists a $k(T)$-isogeny 
(resp. surjective homomorphism, resp. non-zero homomorphism, resp. homomorphism with $\kappa$-dimensional kernel)
$A_T\to B_T$, {{}where  $A_T,$ $B_T$ denote the generic fibres of the base changed abelian schemes $\AAA_T\to T,$ $\BBB_T\to T,$  respectively (cf. Section 3).}
\end{enumerate}
\item[(b)] The following are equivalent:
\begin{enumerate}
\item[(i)] There exists a $\ol{k(S)}$-isogeny (resp. surjective homomorphism, resp. non-zero homomorphism, resp. homomorphism with $\kappa$-dimensional kernel)
$A_{\ol{k(S)}}\to B_{\ol{k(S)}}$.
\item[(ii)] For every $m$-dimensional smooth connected subscheme $T$ of $U$ there exists a $\ol{k(T)}$-isogeny 
(resp. surjective homomorphism, resp. non-zero homomorphism, resp. homomorphism with $\kappa$-dimensional kernel)
$A_{T, \ol{k(T)}}\to B_{T, \ol{k(T)}}$. 
\end{enumerate}
\end{enumerate}
\end{intthm}

\noindent
Ingredients of the proof of Theorem {{}A} include standard methods based on the Tate conjecture (proven by Tate, Zarhin and Faltings cf. Theorem \ref{zf}) 
and some consequences of the Hilbert irreducibility theorem (cf. Lemma \ref{connlemm}), which were inspired by 
Drinfeld's ``conventional formulation of Hilbertianity'' in \cite[Section A.1.]{Drin} {{} and by Section 2 of a recent paper of Cadoret and Tamagawa \cite{CT20}.}  As a formal consequence we obtain:

\begin{intcoro}[Cor. \ref{hilbertmain2}] \label{thmB} Let $S$ be a smooth variety over a finitely generated 
field $k$. Let $\AAA$ be an abelian scheme over $S$ with
generic fiber $A$. Let $U$ be a dense open subscheme of $S$. 
Let $m\in \{0,1,\cdots, \dim(S)\}$.  Assume that $k$ is infinite or that $m\ge 1$. 
\begin{enumerate}
\item[(a)] The following are equivalent:
\begin{enumerate}
\item[(i)] $A$ is not a simple $k(S)$-variety.
\item[(ii)] For every $m$-dimensional smooth connected subscheme $T$ of $U$ 
the fibre $A_T$ is not a simple $k(T)$-variety.
\end{enumerate}
\item[(b)] The following are equivalent:
\begin{enumerate}
\item[(i)] $A_{\ol{k(S)}}$ is not a simple $\ol{k(S)}$-variety.
\item[(ii)] For every $m$-dimensional smooth connected subscheme $T$ of $U$ 
the fibre $A_{T, \ol{k(T)}}$ is not a simple $\ol{k(T)}$-variety
\end{enumerate}
\end{enumerate}
\end{intcoro}

\noindent
Our second main result is the following local to global principle for quadra\-tic isogeny twists of abelian varieties. 
 We 
 call 
 an abelian variety $B/k$ {\sl a quadratic isogeny twist} of an abelian variety $A/k,$ if there exists a quadratic 
twist $A'/k$ of $A$ and a $K$-isogeny $B\to A'$ {{}(cf. Section 2)}. 

\begin{intthm}[Thm. \ref{hilbertmain3}] \label{thmC} Let $S$ be a smooth variety over a finitely generated 
field $k$. Let $\AAA, \BBB$ be abelian schemes over $S$ with
generic fibers $A$ and $B$ respectively. Let $U$ be a dense open subscheme of $S$. 
Let $m\in \{0,1,\cdots, \dim(S)\}$.  Assume that $k$ is infinite or that $m\ge 1$. 
The following are equivalent:
\begin{enumerate}
\item[(a)] $A$ is a quadratic isogeny twist of $B$
\item[(b)] For every $m$-dimensional smooth connected subscheme $T$ of $U$ 
the abelian variety $A_T$ is a quadratic isogeny twist of $B_T$.
\end{enumerate}
The implication (a)$\Rightarrow$(b) holds true also in the case where $k$ is finite and $m=0$. 
\end{intthm}

{{} We remark, that results in Section 4 of the paper are a bit more general than Theorem \ref{thmA}, Corollary \ref{thmB} and Theorem \ref{thmC} 
in that they also cover the situation where $S$ is an arithmetic scheme, but we do not go into the details within this introduction.}

It is clear that in the above statements the case when $k$ is a finite field and $m{=}0$ can not be covered by the Hilbertianity methods. It 
constitutes a separate question which we address under two sets of additional assumptions. 
If  $A$ and $B$ do not 
have nontrivial endomorphisms  geometrically, then we establish a global function field analogue (cf. Proposition 
\ref{Fite-lemm} below)
of Fite's result \cite[Cor. 2.7]{Fite21} following its 
proof quite closely.  By combining Proposition \ref{Fite-lemm} with the 
Hilbertianity approach we augment Theorem \ref{thmC} by the case $k$ finite and $m{=}0.$
 \begin{intthm} [Cor. \ref{Fite-coro}] Let $S$ be a smooth variety over a finite
field $k$. Let $\AAA, \BBB$ be abelian schemes over $S$ with
generic fibers $A$ and $B$ respectively. Let $U$ be a dense open subscheme of $S$. 
Assume that $End_{\ol{k(S)}}(A)=End_{\ol{k(S)}}(B)=\Zz$. 
The following are equivalent:
\begin{enumerate}
\item[(a)]  $A$ is a quadratic isogeny twist of $B$.
\item[(b)] For every closed point $u$ of $U$ the abelian variety
$A_u$ is a quadratic isogeny twist of $B_u$. 
\end{enumerate}
\end{intthm}

\noindent
If abelian varieties $A$ and $B$ meet the so-called {\sl minuscule weights conjecture} of Zarhin (cf. condition MWC, Definition \ref{MWC}), then we apply a 
global function field analogue (cf. Proposition \ref{KL-lemm}) of a result of Khare and Larsen \cite[Thm. 1]{KL20} and prove the following result. {{} It 
completes part (b) of Theorem  \ref{thmA} in case when $k$ is finite and $m{=}0.$} We discuss the current status of Zarhin's 
conjecture in Remark \ref{RMWC}. In particular, it holds true for ordinary abelian varieties over global fields of positive characteristics.   
 
\begin{intthm}[Thm. \ref{KL-thm2}] Let $S$ be a smooth variety over a finite
field $k$. Let $\AAA, \BBB$ be abelian schemes over $S$ with
generic fibers $A$ and $B$ respectively. 
Assume that $A$ satisfies 
$MWC(A)$ and $B$ satisfies $MWC(B)$. 
The following are equivalent:
\begin{enumerate}
\item[(a)] There exists a surjective $\ol{k(S)}$-homomorphism (resp. $\ol{k(S)}$-isogeny) $A_{\ol{k(S)}}\to B_{\ol{k(S)}}$.
\item[(b)] For every closed point $s\in S$ there exists a  surjective  $\ol{k(s)}$-homomor\-phism (resp. $\ol{k(s)}$-isogeny)
$A_{s, \ol{k(s)}}\to B_{s, \ol{k(s)}}$. 
\end{enumerate}
\end{intthm}
\medskip

\noindent
{\bf Structure of the paper.} In Sections 2 and 3 we gathered material which is needed in the sequel including basic facts on: twists of abelian varieties, Galois representations and abelian 
schemes. Section 4 is a central part of the paper. It contains proofs of main results by Hilbertianity methods in the case when $k$ is an infinite field or $m\geq 1.$  In the final section we discuss the remaining case of $k$ finite, $m{=}0$ and work under extra assumptions, either the minuscule weights conjecture or trivial endomorphisms for generic fibres.     

\medskip

\noindent
{\bf Acknowledgments.} The authors were supported by a research grant UMO-2018/31/B/ST1/01474 of the National Centre of Sciences of Poland. S.P. thanks the Mathematics Department at Adam Mickiewicz University in Pozna\'n for hospitality during research visits. We thank J\c edrzej Garnek, Marc Hindry and Bartosz Naskr\c ecki for useful discussions on the topic 
of this paper. 
Finally we want to thank the anonymous referee for a careful reading of the manuscript and for valualble comments.

\newcommand{\Twist}{\mathrm{Twist}}
\newcommand{\Spec}{\mathrm{Spec}}
\section{Preliminaries}
\subsection*{Notation}
For a field $K$ we denote by $\ol{K}$ a separable closure of $K.$ 
If $E/K$ is a Galois extension, we denote by $\Gal(E/K)$ its Galois group and define 
$\Gal(K):=\Gal(\ol{K}/K)$.  A $K$-variety  is a separated algebraic $K$-scheme which is reduced and irreducible. A $K$-curve is a $K$-variety of dimension $1$. For a scheme $S$ and $s\in S$ 
we denote by $k(s)$ the residue field of $s$. Let $n\in \Zz$. Then, as usual, we denote by $S[n^{-1}]$ the
open subscheme of $S$ with underlying set $\{s\in S: n\in k(s)^\times\}$ (where $n$ is viewed as an element 
of $k(s)$ via the ring homomorphism $\Zz\to k(s)$). We let $\Ll$ be the set of all rational primes and define $\Ll(S):=\{\ell\in\Ll: S[\ell^{-1}]\neq \emptyset\}$. 
If $S$ is reduced and irreducible we denote, following EGA, by $R(S)$ the
function field of $S$. If $S$ is a $K$-variety, we sometimes write $K(S)$ instead of $R(S)$. 
If $T$ is a finite free $\Zz_\ell$-module, $V=T\otimes_{\Zz_\ell} \Qq_\ell$ and $\Gamma$ is a subgroup $\GL_{T}(\Zz_\ell)$, then we denote by $\Gamma^{\Zar}$
the Zariski closure of $\Gamma$ inside the algebraic group $\GL_{V}/\Qq_\ell$ so that $\Gamma^{\Zar}$ is an algebraic group of $\Qq_\ell$.  If $\ul{G}$ is an algebraic group over $\Qq_\ell$, then we denote by
$\ul{G}^\circ$ the connected component of the identity element of $\ul{G}$. 

\subsection*{Twists of abelian varieties.}  
Let $K$ be a field and $A$ and $B$ abelian varieties over $K$. Let $E/K$ be a Galois extension. 
 We call $B$ an $E/K$-twist of $A$ if there is an $E$-isomorphism $A_E\to B_E$.   In subsequent sections we are mainly interested in the case 
 $E=\ol{K}$.
 We denote by 
$\Twist_{E/K}(A)$ the set of all isomorphism classes of $E/K$-twists of $A$ and define $\Twist(A):=\Twist_{\ol{K}/K}(A)$. 
There are natural operations $\rho_A: \Gal(K)\to \Aut_K(A_E)$ and $\rho_B: \Gal(K)\to \Aut_K(B_E)$ and an operation
$$\Gal(K)\times \Hom_E(A_E, B_E)\to \Hom_E(A_E, B_E),\ {}^\sigma f:=\rho_B(\sigma)\circ f\circ \rho_A(\sigma){{}^{-1}}.$$
Now assume that $B/K$ is an $E/K$-twist of $A$ and choose an $E$-isomorphism $f: A_E\to B_E$. Then 
$$\xi: \Gal(E/K)\to \Aut_E(A_E),\ \xi(\sigma)=f^{-1}\circ {}^\sigma f$$
is a $1$-cocycle (i.e. $\xi(\sigma\eta)=\xi(\sigma)\circ{}^\sigma\xi(\eta)$) whose cohomology class does not depend on the choice of $f$. 
Let $$\tilde{\rho}_B: \Gal(E/K) \to \Aut_K(A_E),\ \tilde{\rho}_B(\sigma)=f^{-1}\circ  \rho_B(\sigma)\circ f$$
be the operation on $A$ derived from $\rho_B$ via transport of structure via $f$. Then an easy calculation shows that
\begin{align}\label{twist:rep1}
\tilde{\rho}_B(\sigma)=\xi(\sigma) \rho_A(\sigma)\quad \forall\sigma\in \Gal(E/K)
\end{align}
is simply the action $\rho_A$ twisted by the $1$-cocycle $\xi$. 
It is well-known that this sets up a 
bijective map
$$\Xi_{E/K}: \Twist_{E/K}(A)\to H^1(\Gal(E/K), \Aut_E(A_E)).$$
There is a natural map
$$j: H^1(\Gal(E/K), \{\pm Id_A\})\to H^1(\Gal(E/K), \Aut_{E}(A_{E})).$$
For a quadratic character $\chi\in \Hom(\Gal(E/K), \{\pm Id_A\})$ we denote by $A_\chi$ an $E/K$-twist of $A$
corresponding to $j(\chi)$ under $\Xi_{E/K}$. We call $B$ {\sl a quadratic isogeny twist of $A,$} if there exists a quadratic character $\chi\in \Hom(\Gal(\ol{K}/K), \{\pm Id_A\})$
and a $K$-isogeny  $B\to A_\chi$.  Sometimes we tacitly identify $\{\pm Id_A\}$ with $\mu_2(\Qq)=\{\pm 1\}$. 

\begin{rema} Let $A$ be an abelian variety over a {\em finite} field $K$. Then there is a unique non-trivial quadratic character $\Gal(K)\to \{\pm Id_A\}$ and accordingly 
a unique non-trivial quadratic twist of $A$. This is clear since $\Gal(K)\cong \hat{\Zz}$.
\end{rema} 

\noindent
\subsection*{Galois representations.} Let $K$ be a field of characteristic $p\ge 0$ and $A/K$ an abelian variety. For a rational prime $\ell\neq p$ we denote by
$$\rho_{A, \ell}: \Gal(K)\to \GL_{V_\ell(A)}(\Qq_\ell)$$
the corresponding $\ell$-adic Galois representation,  where $T_\ell(A)$ is the $\ell$-adic Tate 
module of $A$ and $V_\ell(A)=T_\ell(A)\otimes_{\Zz_\ell} \Qq_\ell$.  
Then the Zariski closure $\rho_{A, \ell} (\Gal(K))^{\Zar}$ of $\rho_{A, \ell}(\Gal(K))$ in $\GL_{V_\ell(A)}$ and its identity component 
$(\rho_{A, \ell} (\Gal(K))^{\Zar})^\circ$ 
are algebraic groups over $\Qq_\ell$. 

\newcommand{\Fr}{\mathrm{Fr}}

If $K$ is {finite}, then we define the $L$-series of $A$ by
$$L(A/K, T):=\det(Id_{V_\ell (A)}-  \rho_{A,\ell}(\Fr)T)$$
where $\Fr\in \Gal(K)$ is the Frobenius element. 
This characteristic polynomial $L(A/K,T)$ has integer coefficients and does not
depend  on the rational prime $\ell\neq p$ by the Weil conjectures. We have the following elementary but useful fact. 

\begin{lemm} \label{isogeny1} Let $K$ be a field of characteristic $p\ge 0$ and $\ell$ a prime different from $p$. Let $A/K$ and $B/K$ be abelian varieties. Let $f: A\to B$ be a homomorphism 
and let $V_\ell(f): 
V_\ell(A)\to V_\ell(B)$ be the homomorphism of $\Qq_\ell$-vector spaces induced by $f$. Then: 
\begin{enumerate}
 \item[a)] $\dim_{\Qq_\ell}(\im(V_\ell (f)))=2\dim(\im(f))$.
 \item[b)] $\dim_{\Qq_\ell}(\ker(V_\ell (f)))=2\dim(\ker(f))$.
 \item[c)]
 $f: A\to B$ is an isogeny if and only if the map $V_\ell(f)$ is bijective. In particular, we have $\rho_{A,\ell}\cong \rho_{B,\ell}$ provided $B$ is $K$-isogenous to $A$. 
 \end{enumerate}
\end{lemm}

{\em Proof.} We note that $T_\ell(A)= \Hom(\Qq_\ell/\Zz_\ell, A(\ol{K}))$. For the purpose of that proof we put 
 $T_\ell(M):=\Hom(\Qq_\ell/\Zz_\ell, M)$ 
for an arbitrariy abelian group $M$. 
Let $C=\ker(f)$ 
Let $I=\im(f)$. Then $C^\circ$ and $I$ are abelian varieties. 
From the exact sequence $$0\to C^\circ(\ol{K})\to C(\ol{K}) \to C(\ol{K})/C^\circ(\ol{K})\to 0$$ we derive an exact sequence
$$0\to T_\ell (C^\circ(\ol{K}))\to T_\ell (C(\ol{K})) \to T_\ell (C(\ol{K})/C^\circ(\ol{K}))$$
because the functor $\Hom(\Qq_\ell/\Zz_\ell, -)$ is left exact. But $T_\ell (C(\ol{K})/C^\circ(\ol{K}))=0$ because the group $C(\ol{K})/C^\circ(\ol{K})$ is finite. Thus
$T_\ell (C^\circ(\ol{K}))= T_\ell (C(\ol{K}))$.  From the exact sequence 
\begin{align}\label{isogeny1:eq1}
0\to C(\ol{K})\to A(\ol{K})\to I(\ol{K})\to 0
\end{align}
we obtain an exact sequence {{}$0\to T_\ell (C(\ol{K}))\to T_\ell (A(\ol{K}))\to T_\ell (I(\ol{K}))$}, again because the functor $\Hom(\Qq_\ell/\Zz_\ell, -)$ is left exact. Using 
$T_\ell (C^\circ(\ol{K}))= T_\ell (C(\ol{K}))$ and tensoring with $\Qq_\ell$ we see that we have an exact sequence 
\begin{align}\label{isogeny1:eq2}
0\to V_\ell (C^\circ) \to V_\ell (A)\to V_\ell (I).
\end{align}
Let $c=\dim(C)$, $a=\dim(A)$ and $i=\dim(I)$. Then $a=i+c$ by \eqref{isogeny1:eq1}. Also $\dim_{\Qq_\ell}(V_\ell(C))=2c$, 
$\dim_{\Qq_\ell}(V_\ell(A))=2a$ and $\dim_{\Qq_\ell}(V_\ell(I))=2i$ (cf. \cite[Remark 8.4]{milneav}). The image of the right hand map of the exact sequence
\eqref{isogeny1:eq2} has dimension $2a-2c=2i$, too, hence $V_\ell(A)\to V_\ell(I)$ is surjective. By the left exactness of $\Hom(\Qq_\ell/\Zz_\ell, -)$ we also see that 
$V_\ell (I)\to V_\ell(B)$ is injective. Thus $V_\ell(I)=\im(V_\ell(f))$ and $V_\ell(C^\circ)=\ker(V_\ell(f))$. It follows that 
\begin{align*}
&\dim_{\Qq_\ell}(\ker(V_\ell(f))){=}\dim_{\Qq_\ell}(V_\ell (C^\circ)){=}2\dim(C^\circ){=2}\dim(C){=}2\dim(\ker(f)),\\
& \dim_{\Qq_\ell}(\im(V_\ell(f))){=}\dim_{\Qq_\ell}(V_\ell (I)){=}2\dim(I){=}2\dim(\im(f)),
\end{align*}
as desired. This finishes the proof of a) and b). For c)  put $b:=\dim(B)$ and note that:
{{}\begin{align*}
\mbox{$f$ is an isogeny}&\Leftrightarrow \mbox{$f$ is surjective with a finite kernel}\\
&\Leftrightarrow \dim(\ker(f)){=}0\,\,{\rm and} \,\dim(\im(f)){=}b\\
&\Leftrightarrow \dim_{\Qq_\ell}(\ker(V_\ell (f))){=}0\,\, {\rm and} \,\dim_{\Qq_\ell}(\im(V_\ell (f))){=}2b\\
& \Leftrightarrow \mbox{$V_\ell (f)$ is bijective.}
\end{align*}}\hfill $\Box$
\medskip

\noindent
There is the following celebrated result due to Faltings, Tate and Zarhin.

\begin{thm}(cf. \cite{Fal83}, \cite{Tate66}, \cite{Zar14}) \label{zf}
Let $K$ be a finitely generated field of characteristic $p\ge 0$ and $A/K$ and $B /K$ be abelian varieties.
\begin{itemize}
\item[(a)] $\rho_{A, \ell}$ is semisimple for every rational prime $\ell\neq p$.
\item[(b)] The natural homomorphism
$$\Hom_K(A,B)\otimes \Zz_\ell\to \Hom_{\Gal(K)}(T_\ell(A), T_\ell(B))$$
is bijective. 
\end{itemize}
\end{thm}

\begin{coro}\label{zfcoro}
Let $K$ be a finitely generated field of characteristic $p\ge 0$ and $A/K$ and $B /K$ be abelian varieties. The following statements are equivalent. 
\begin{enumerate}
\item[(a)] There exists a $K$-isogeny $f: A\to B$.
\item[(b)] There exists an isomorphism $\rho_{A, \ell}\cong \rho_{B, \ell}$ of $\Qq_\ell$-representations of $\Gal(K)$. 
\end{enumerate}
If $K$ is finite, then conditions (a) and (b) are also equivalent to the following condition
\begin{enumerate}
\item[(c)] $L(A/K, T)=L(B/K, T)$. 
\end{enumerate}
\end{coro}

{\em Proof.} The implication (a)$\Rightarrow$(b) is Lemma \ref{isogeny1}. The implication (b)$\Rightarrow$(a) is a consequence of Theorem \ref{zf}, cf. \cite[Korollar 2]{Fal83}. If $K$ is finite, then the equivalence
(b)$\Leftrightarrow$(c) has been established by Tate \cite[Theorem 1]{Tate66}.\hfill $\Box$

\begin{rema} \label{twistrep1} Let $K$ be a finitely generated field of characteristic $p\ge 0$ and $A/K$ an abelian variety. Let $\chi: \Gal(K)\to \{\pm Id_A\}$ be
a quadratic character. Then:
 \begin{enumerate}
\item[(a)] $\rho_{A_\chi, \ell}\cong \chi \otimes \rho_{A, \ell}$,
\item[(b)] If $K$ is finite and $\chi$ is the (unique) non-trivial  quadratic character of $\Gal(K)$, then $L(A_\chi/K,T)=L(A/K, -T)$.  
\end{enumerate}
\end{rema}

{\em Proof.} (a) is immediate from equation \eqref{twist:rep1}. Part (b) is an immediate consequence thereof. \hfill $\Box$

\begin{lemm} \label{TwistLemm} Let $K$ be a finitely generated field of characteristic $p\ge 0$. 
Let $A,B$ be abelian varieties over $K$ and $\ell\neq p$ a rational prime. Then the following statements are equivalent.
\begin{enumerate}
\item[(a)] $B$ is a quadratic isogeny twist of $A$.
\item[(b)] There exists a quadratic character $\chi: \Gal(K)\to \{\pm Id_A\}$ such that $\rho_{B, \ell}\cong \chi\otimes \rho_{A, \ell}$. 
\end{enumerate}
If $K$ is finite, then the equivalent conditions (a) and (b) are also equivalent to
\begin{enumerate}
\item[(c)] $L(B/K, T)=L(A/K, T)$ or $L(B/K,T)= L(A/K, -T)$. 
\end{enumerate}
\end{lemm} 

{\em Proof.} We prove the implication (a)$\Rightarrow$(b): Assume $B$ is a quadratic isogeny twist of $A$. Then there exists a quadratic character
$\chi: \Gal(K)\to \{\pm Id_A\}$ and a $K$-isogeny $B\to A_\chi$. From Lemma \ref{isogeny1} and Remark \ref{twistrep1} we conclude that
$$\rho_{B, \ell}\cong\rho_{A_\chi, \ell}\cong \chi\otimes \rho_{A, \ell}.$$

For the proof of  the implication (b)$\Rightarrow$(a) assume that there 
exists a quadratic character $\chi: \Gal(K)\to \{\pm Id_A\}$ such that $\rho_{B, \ell}\cong \chi\otimes \rho_{A, \ell}$. From Remark \ref{twistrep1} we conclude
that $\rho_{B, \ell}\cong \rho_{A_\chi, \ell}$ and then Corollary \ref{zfcoro} implies that there exists an isogeny $B\to A_\chi$. It follows that $B$ is a quadratic
isogeny twist of $A$.

From now on assume that $K$ is finite. Then $\Gal(K)=\hat{\Zz}$ and thus there exist exactly two quadratic characters thereof. 

We prove the implication (b)$\Rightarrow$(c): If $\chi$ is the trivial character then $\rho_{B, \ell}\cong \rho_{A, \ell}$ and hence $L(A/K, T)=L(B/K,T)$. 
If $\chi$ is non-trivial, then $\rho_{B, \ell}\cong \rho_{A_\chi, \ell}$, hence
$$L(B/K, T)=L(A_\chi, T)=L(A/K,-T)$$
by Corollary \ref{zfcoro} and Remark \ref{twistrep1}.

 We prove the implication (c)$\Rightarrow$(a): If $L(B/K,T)=L(A/K,T)$, then $A$ is isogenuous to $B$ by Corollary \ref{zfcoro}. If 
 $L(B/K,T)=L(A/K,-T)$ and $\chi$ the nontrivial quadratic character, then $L(B/K,T)=L(A_\chi/K,T)$ and Corollary \ref{zfcoro} implies that
 $B$ is isogenous to $A_\chi$, as desired. \hfill $\Box$
 
 \begin{rema} \label{rema:tatehom} Let 
 $A$ and $B$ be abelian varieties over a 
 field $K$ and $\ell$ a prime different from the characteristic of $K$. Note that $\rho_{A\times B, \ell}(\sigma)=\rho_{A,\ell}(\sigma)\times
 \rho_{B,\ell}(\sigma)$ for all $\sigma\in\Gal(K)$. 
  Thus 
 {$$\rho_{A\times B, \ell}(\Gal(K))\subset \GL_{V_\ell (A)}(\Qq_\ell)\times \GL_{V_\ell (B)}(\Qq_\ell)\subset 
 \GL_{V_\ell (A)\times V_\ell (B)}(\Qq_\ell).$$
 Let $p_A$ and $p_B$ be the projections of the product $\GL_{V_\ell (A)}(\Qq_\ell)\times \GL_{V_\ell (B)}(\Qq_\ell)$. 
 Then $\rho_{A, \ell}(\sigma)=p_A(\rho_{A\times B, \ell}(\sigma))$ and $\rho_{B, \ell}(\sigma)=p_B(\rho_{A\times B, \ell}(\sigma))$ 
 for all $\sigma \in\Gal(K)$. Thus the actions of $\Gal(K)$  on $V_\ell(A)$ and $V_\ell(B)$ and the induced action on 
 $\Hom_{\Qq_\ell}(V_\ell(A), V_\ell(B))$ factor through $G=\rho_{A\times B, \ell}(\Gal(K))$.
 Let $\uG=\rho_{A\times B, \ell}(\Gal(K))^{\Zar}$.  If $K$ is finitely generated, then the natural map
 \begin{align}\label{eq:longHom}
 \Hom_K(A,B)\otimes_{\Zz} \Qq_\ell\to \Hom_{\Qq_\ell}(V_\ell(A), V_\ell(B))^{{G}}= 
 \Hom_{\Qq_\ell}(V_\ell (A), V_\ell (B))^{\uG}
 \end{align}
 is bijective by the Tate conjecture (cf. Theorem \ref{zf}).}
\end{rema}

The following lemma and its proof are inspired by \cite[Prop 2.10]{SZ95}. 
\begin{lemm}\label{lemm:tatevariant}
Let $A$ and $B$ be abelian varieties over a field $K$ and keep the notation from Remark \ref{rema:tatehom}.
Let $K^\circ$ be the fixed field of the group $(\rho_{A\times B,\ell})^{-1}(\uG^\circ)$. Then 
\begin{align}\label{eq:tatevariant}
\Hom_{K^\circ}(A_{K^\circ}, B_{K^\circ})= \Hom_{\ol{K}}(A_{\ol{K}}, B_{\ol{K}}).
\end{align} If $K$ is finitely generated, then the 
above specialization map
\eqref{eq:longHom} induces a bijective map
 \begin{align}\label{eq:tateinvariant2}
 \Hom_{\ol{K}}(A_{\ol{K}},B_{\ol{K}})\otimes_\Zz \Qq_\ell \to \Hom_{\Qq_\ell}(V_\ell (A), V_\ell (B))^{\uG^\circ}.
 \end{align}
\end{lemm}

{\em Proof.} For \eqref{eq:tatevariant} it  is enough to prove that $$\Hom_{K^\circ}(A_{K^\circ}, B_{K^\circ})= \Hom_{K'}(A_{K'}, B_{K'})$$ for every finite Galois extension 
$K'/K$ with $K^\circ \subset K'$. 
The profinite group $\Gal(K')$ is of finite index in $\Gal(K)$ and thus $\rho_{A\times B,\ell}(\Gal(K'))^{\Zar}$ is a closed subgroup scheme of finite index in $\uG^\circ$. Thus finitely many cosets of 
$\rho_{A\times B,\ell}(\Gal(K'))^{\Zar}$ cover $\uG^\circ$. As $\uG^\circ$ is connected, it follows that 
 $$\rho_{A\times B,\ell}(\Gal(K'))^{\Zar}=\uG^\circ.$$

Let $H= \Hom_{\Qq_\ell}(V_\ell (A), V_\ell (B))$ and consider the diagram
 $$\xymatrix{
\Hom_{K^\circ}(A_{K^\circ},B_{K^\circ})\ar[r]\ar[d] &H^{\uG^\circ}\ar[d]_{Id}\ar@{=}[r]& H^{\Gal(K^\circ)}\ar[d]_{Id}\\
\Hom_{K'}(A_{K'},B_{K'}) \ar[r] & H^{\uG^\circ}\ar@{=}[r]&H^{\Gal(K')}}$$
The horizontal arrows are injective 
(cf. Theorem \cite[Prop. 12.2]{milneav}). Thus, for every $f\in \Hom_{K'}(A_{K'},B_{K'})$, the map $V_\ell(f)\in H$ is invariant unter $\Gal(K^\circ)$ and, by the
injectivity of the horizontal maps of the diagram, $f$ itself is invariant under $\Gal(K^\circ)$ which impies $f\in \Hom_{K^\circ}(A_{K^\circ},B_{K^\circ})$. 
This finishes up the proof of \eqref{eq:tatevariant}. If $K$ is finitely generated, then \eqref{eq:tatevariant} together with the Tate conjecture 
(cf. {Theorem \ref{zf}}) implies that 
\eqref{eq:tateinvariant2} is bijective. \hfill{$\Box$}

 \section{Abelian schemes}
 
 {\em Throughout this section let $S$ and $T$ be noetherian,  normal and connected schemes. 
 Note that then the local rings of $S$ are domains and thus \cite[6.1.10]{EGAI} (together with the connectedness of $S$) implies that $S$ is in fact irreducible.  Furthermore $S$ is reduced. Similarly $T$ 
 is reduced and irreducible.
 Let $F$ (resp. $E$) be
 the function field of $S$ (resp. $T$). Let $u: T\to S$ be a morphism. Furthermore consider the point 
 $$t: \Spec(E)\to T\buildrel u\over \longrightarrow S$$
 of $S$.  We fix throughout this section a rational prime $\ell\in \Ll(T)$. This is possible by the following 
 lemma. We note that automatically $\ell\in \Ll(S)$.}
 
 \begin{lemm}\label{ellex}
For every scheme $X,$ the set $\Ll(X)$ is non-empty.
 \end{lemm}
 
 {\em Proof.} If there exists $p\in \Ll\cup \{0\}$ such that $\chara(k(x))=p$ for all $x\in X$, then $\Ll(X)=\Ll
 \setminus \{p\}$ is not empty.  Otherwise,  there exist $x_1, x_2\in X$ such that $\chara(k(x_1))\neq \chara(k(x_2))$
 and $\chara(k(x_1))>0$, hence putting $\ell:=\chara(k(x_1))$ we get $x_2\in X[\ell^{-1}]$ and thus $\ell\in \Ll(X),$ so that $\Ll(X)$ is non-empty also in that case. \hfill $\Box$ 
 \par\medskip
 
Let $\AAA/S$ be 
an abelian scheme with generic fibre $A/F$. Let $\AAA_T:=\AAA\times_S T$ and let $A_T:=\AAA\times_S \Spec(E)$ be the generic fibre
of $\AAA_T\to T$. We then have cartesian squares

$$\xymatrix{A_T \ar[r]\ar[d]&  \AAA_T\ar[r]\ar[d] & \AAA\ar[d] & A\ar[l]\ar[d]\\
\Spec(E)\ar[r]&  T\ar[r]^u &S & \Spec(F)\ar[l]\\
}$$

Let $n\in\Nn$.
Then the restriction of $\AAA[\ell^n]$ to $S[\ell^{-1}]$ is a finite \'etale $S[\ell^{-1}]$-scheme of rank $\ell^{2\dim(A)n}$. In particular the action of $\Gal(F)$ on $A[\ell^n](\ol{F})$ factors through
$\pi_1(S[\ell^{-1}])$. Likewise $\rho_{A, \ell}$ factors through $\pi_1(S[\ell ^{-1}])$, i.e. we 
can view it as a homomorphism $\rho_{A, \ell}: \pi_1(S[\ell^{-1}])\to GL_{V_\ell(A)}(\Qq_\ell)$. 
There exists a finite connected \'etale cover $S'\to S[\ell^{-1}]$ such that $\AAA[\ell^n]\times S'$ splits up in a coproduct of $\ell^{2n\dim(A)}$ copies of $S'$,  i.e. ,  it is a constant group scheme. There is a point $t': \Spec(\ol{E})\to S'$ over $t$ and  a point $\xi': \ol{F}\to S'$ lying over the generic point of $S$. As $\AAA[\ell^n]\times_S S'$ is a constant group scheme,  the natural maps 
 $$ A[\ell^n](\ol{F})\leftarrow \AAA[\ell^n](S')\rightarrow A_T[\ell^n](\ol{E})$$
 derived from these points are bijective. We thus get a specialization isomorphism
 $$s_{A,T,\ell^n}: A[\ell^n](\ol{F})\to A_T[\ell^n](\ol{E})$$
 and accordingly a specialization isomorphism
 $$s_{A,T,\ell^\infty}: T_\ell(A)\to T_\ell(A_T).$$ 
 These are equivariant for the action of $\pi_1(S[\ell^{-1}])$.  In particular we have $\dim(A)=\dim(A_T)$.

 \begin{defi} Let $G$ be a group, $S_0$ an open subscheme of $S$ and $\rho: \pi_1(S_0)\to G$ 
 a homomorphism. Let $T_0:=u^{-1}(S_0)$. We define $u^*\rho$ to be the composite 
 morphism\footnote{Note: The fundamental group of the empty scheme is the trivial group.} 
 $$\pi_1(T_0)\buildrel u_*\over \longrightarrow \pi_1(S_0)\to G.$$
If $u$ is clear from the context we put $\rho_T:=u^*\rho$. We say that $T$ (or $u$) is {\bf $\rho$-generic}
if $T_0$ is not empty and $\rho_T(\pi_1(T_0))=\rho(\pi_1(S_0))$.
\end{defi}

We often apply this notation in the case where $S_0=S[\ell^{-1}]$ and $\rho=\rho_{A,\ell}$. Note that in
that case $u^{-1}(S_0)=T[\ell^{-1}]$ is automatically non-empty by our choice of $\ell$. 

\begin{rema}
The specialization isomorphism $s_{A,T,\ell^\infty}: T_\ell(A)\to T_\ell(A_T)$ gives an isomorphism 
$$u^* \rho_{A,\ell}\cong \rho_{A_T, \ell}.$$
of representations of $\Gal(E)$. 
\end{rema}
 
{\em For the rest of this section let $\BBB$ be an abelian scheme over $S$ with generic fibre $B/F$. Define $\BBB_T:=\BBB\times_S T$ and 
let $B_T:=\BBB\times_T \Spec(E)$ be the generic fibre of $\BBB_T\to T$. }

\begin{lemm} \label{CF-lemm}  The canonical 
map 
$$\Hom_S(\AAA,\BBB)\to \Hom_F(A,B),\ f\mapsto f_F$$
is bijective. 
\end{lemm}

{\em Proof.}  Let $\UUU$ be the category of all dense affine open subschemes of $S$. Then $\Spec(F)=\plim_{U\in\UUU} U$ (projective limit of schemes, cf. \cite[8.2]{EGAIV3}). It follows from
\cite[8.8.2]{EGAIV3} that the natural map
\begin{align}\label{CF-bij0}
\ilim_{U\in \UUU}\Hom_U(\AAA_U, \BBB_U)\to \Hom_F(A,B),
\end{align}
induced by the natural maps $\Hom_U(\AAA_U, \BBB_U)\to \Hom_F(A,B),\ f\mapsto f_F$, is bijective.
Moreover, by \cite[Prop. 2.7]{CF}, for every  $U\in \UUU$ the natural map 
\begin{align*}
\Hom_S(\AAA,\BBB)\to \Hom_U(\AAA_U, \BBB_U),\ f\mapsto f_U
\end{align*}
is bijective, and therefore the natural map
\begin{align}\label{CF-bij2}
\Hom_S(\AAA,\BBB)\to \ilim_{U\in\UUU} \Hom_U(\AAA_U, \BBB_U)
\end{align}
is bijective.
The assertion is immediate from the bijectivity of \eqref{CF-bij0} and \eqref{CF-bij2}. 
\hfill $\Box$
 
 \begin{lemm} \label{rig-lemm} 
 The canonical map
 \begin{align}\label{reszero}
 \Hom_S(\AAA, \BBB)\to \Hom_T(\AAA_T, \BBB_T)
 \end{align}
 is injective.
 \end{lemm}

 {\em Proof.} Choose a point $t\in T$. 
 We even prove that the composite homomorphism 
 \begin{align}\label{reszero2}\Hom_S(\AAA, \BBB)\to \Hom_T(\AAA_T, \BBB_T)\to \Hom_{k(t)}(\AAA_t, \BBB_t)
 \end{align}
 is injective.
 Let $f$ be in the kernel of $\eqref{reszero2}$. Then $$f_t: \AAA\times_S \Spec(k(t))\to \BBB\times_S \Spec(k(t))$$
 is the zero homomorphism. Hence $f$ must be the zero homomorphism by the rigidity 
 lemma \cite[20.1]{milneav}.
 \hfill $\Box$
 \medskip

 Note that, by the above Lemmata \ref{CF-lemm} and \ref{rig-lemm}, there is a canonical injective specialization map
\begin{align}\label{spz1}
r_{\AAA,\BBB,T, S}: \Hom_{F}(A,B){\cong}\Hom_S(\AAA,\BBB)\longrightarrow\Hom_T(\AAA_T, \BBB_T){\cong}& \\
\notag \cong 
 \Hom_{E}(A_T, B_T).
\end{align}
Let 
$\rho_\ell:=\rho_{A\times B,\ell}$, $\uG:=\rho_\ell(\pi_1(S[\ell^{-1}]))^{\Zar}$ and let $S'$ be 
the finite \'etale cover of $S[\ell^{-1}]$ corresponding to the subgroup $\rho_\ell^{-1}(\uG^\circ)$ of $\pi_1(S[\ell^{-1}])$. Let $F'$ be the function field of $S'$. 
Then, by Lemma \ref{lemm:tatevariant}, 
$\Hom_{F'}(A_{F'}, B_{F'})=\Hom_{\ol{F}}(A_{\ol{F}},B_{\ol{F}})$. 
Let $T'$ be an irreducible component of $T\times_S S'$ and $E'$ the function field of $T'$. Then $E'/E$ is a finite separable extension and we can consider the composite map 
\begin{align}\label{rbar}
\ol{r}_{\AAA,\BBB,T, S}: &\Hom_{\ol{F}}(A_{\ol{F}},B_{\ol{F}})=\Hom_{F'}(A_{F'}, B_{F'})
{{}\buildrel r_{\AAA,\BBB,T',S'} \over\longrightarrow}\\\notag
&\to \Hom_{E'}(A_{T,E'}, B_{T,E'})\to \Hom_{\ol{E}}(A_{T,\ol{E}}, B_{T,\ol{E}})
\end{align}
which tacitly depends on the choice of $T'$ and on the choice of an embedding  $E'\to \ol{E}$.
This map $$\ol{r}_{\AAA,\BBB,T, S}: \Hom_{\ol{F}}(A_{\ol{F}},B_{\ol{F}})\to \Hom_{\ol{E}}(A_{T,\ol{E}}, B_{T,\ol{E}})$$ is injective.

\begin{lemm} \label{specializehoms} Let $f: A\to B$ be a homomorphism. Let $f_T: A_T\to B_T$ be the homomorphism $f_T:=r_{\AAA,\BBB,T,S}(f)$ 
obtained from $f$ by specialization. 
\begin{enumerate}
\item[(a)] $\dim(\im(f))=\dim(\im(f_T))$,
\item[(b)] $\dim(\ker(f))=\dim(\ker(f_T))$,
\item[(c)] $f$ has finite kernel (resp. is surjective, resp. is an isogeny) if and only if $f_T$ has finite kernel (resp. is surjective, resp. is an isogeny).
\end{enumerate}
\end{lemm}
 
 {\em Proof.} 
 With the above specialization isomorphisms we construct a diagram
 $$\xymatrix{V_\ell(A) \ar[r]^{V_\ell(f)} \ar[d] &  V_\ell(B) \ar[d]\\
V_\ell(A_T) \ar[r]^{V_\ell(f_T)}& V_\ell(B_T)
}$$
 whose vertical arrows are bijective. Together with Lemma \ref{isogeny1} we get:
 \begin{align*}
& \dim(\ker(f))=\frac{1}{2}\dim(\ker(V_\ell(f)))=\frac{1}{2}\dim(\ker(V_\ell(f_T)))= \dim(\ker(f_T)),\\
 &\dim(\im(f))=\frac{1}{2}\dim(\im(V_\ell(f)))=\frac{1}{2}\dim(\im(V_\ell(f_T)))= \dim(\im(f_T)).
 \end{align*}
 This proves (a) and (b), and (c) is immediate from that.  \hfill $\Box$. 
 
 \begin{rema} \label{rema:spec} For $f\in \Hom_{\ol{F}}(A_{\ol{F}},B_{\ol{F}})$ and $f_T:=\ol{r}_{\AAA,\BBB,T,S}(f)$ one can compare dimension data of
 $f$ and $f_T$ in a completely analogous way. This is a consequence of Lemma \ref{specializehoms} and the construction of
 $\ol{r}_{\AAA,\BBB,T,S}.$
 \end{rema}

\begin{lemm}\label{fingenhomspez} 
If both fields $E$ and $F$ are finitely generated and $u: T\to S$ is $\rho_{A\times B, \ell}$-generic, 
then the canonical maps
$$r_{\AAA,\BBB,T,S}\otimes \Zz_\ell: \Hom_{F}(A,B)\otimes \Zz_\ell \to \Hom_{E}(A_T, B_T)\otimes \Zz_\ell$$
and 
$$\ol{r}_{\AAA,\BBB,T,S}\otimes \Zz_\ell: \Hom_{\ol{F}}(A_{\ol{F}},B_{\ol{F}})\otimes \Zz_\ell \to 
\Hom_{\ol{E}}(A_{T,\ol{E}}, B_{T,\ol{E}})\otimes \Zz_\ell$$
are bijective.  In particular,  under these assumptions $\coker(r_{\AAA,\BBB,T,S})$ and $\coker(\ol{r}_{\AAA,\BBB,T,S})$  are finite groups of order prime to $\ell.$
\end{lemm}

{\em Proof.} We identify, $T_\ell (A)$ with $T_\ell (A_T)$ and $T_\ell( B)$ with $T_\ell (B_T)$ along the natural and equivariant specialization 
isomorphisms.  Let $\rho_\ell:=\rho_{A\times B, \ell}$, 
$\ul{G}:=\rho_\ell(\pi_1(S[\ell^{-1}]))^{\Zar}$ and $\ul{G}_T:=(\rho_\ell)_T(\pi_1(T[\ell^{-1}]))^{\Zar}$. 
Let $$r_{\ell}=r_{\AAA,\BBB,T,S}\otimes \Zz_\ell$$ and $\ol{r_{\ell}}=\ol{r}_{\AAA,\BBB,T,S}\otimes \Zz_\ell$. 
We then have a commutative diagram
$$\xymatrix{\Hom_{F}(A,B)\otimes {\Zz_\ell}\ar[r]^{r_{\ell}} \ar[d] &  \Hom_{E}(A_T,B_T)\otimes {\Zz_\ell}\ar[d]\\
\Hom_{\Zz_\ell}(V_\ell (A), V_\ell (B))^{\ul{G}} \ar[r]& \Hom_{\Zz_\ell}(V_\ell (A), V_\ell (B))^{\ul{G}_T},
}$$
where the vertical maps are bijective (cf. Tate conjecture,  Remark \ref{rema:tatehom}). 
The lower horizonal map is bijective because
$\ul{G}=\ul{G}_T$ by our assumption that $T$ is $\rho_\ell$-generic. Hence the upper horizontal map $r_{\ell}$
is bijective,  too.  Next consider the commutative diagram
$$\xymatrix{\Hom_{\ol{F}}(A_{\ol{F}},B_{\ol{F}})\otimes {\Zz_\ell}\ar[r]^{\ol{r_{\ell}}} \ar[d] &  \Hom_{\ol{E}}(A_{T,\ol{E}},B_{T,\ol{E}})\otimes {\Zz_\ell}\ar[d]\\
\Hom_{\Zz_\ell}(V_\ell (A), V_\ell (B))^{\ul{G}^\circ} \ar[r]& \Hom_{\Zz_\ell}(V_\ell (A), V_\ell (B))^{\ul{G}_T^\circ}
}$$
where the vertical maps are bijective (cf. Tate conjecture,  Lemma \ref{lemm:tatevariant}). 
The lower horizonal map is bijective because
$\ul{G}^\circ=\ul{G}_T^\circ$. Hence, $\ol{r_{\ell}}$ is bijective. The $\Zz$-module $\Hom_{\ol{E}}(A_{T,\ol{E}}, B_{T,\ol{E}})$ is free and finitely generated. 
Hen\-ce, the statement about the cokernels 
of $r_{\AAA,\BBB,T,S}$ and $\ol{r}_{\AAA,\BBB,T,S}$ follows as well.\hfill $\Box$

\section{Hilbertianity}
 {\em Throughout this section let $\Bb$ be a regular noetherian connected scheme. Let $S$ be a connected scheme and $f: S\to \Bb$ a 
 morphism of finite type which is assumed to be smooth.  Note that 
 $S$ and $Z$ are reduced and irreducible.
 Let $d$ be the relative dimension of $S/\Bb$. Let $K$ be the function field of $\Bb$. Assume that $K$ is finitely generated. 
We denote by $\Sm_m(S/\Bb)$ the set of all connected subschemes $T$ of $S$ such that the restriction 
 $f|T:T\to \Bb$ is smooth of relative dimension $m$. Note that every $T\in \Sm_m(S/\Bb)$ is regular and connected, hence reduced and irrecucible.}

The aim of this section is to consider specializations to subvarieties in $\Sm_m(S/\Bb)$. 
The results of this section are in our opinion most interesting in the following cases:
\begin{enumerate}
\item The case where $\Bb=\Spec(k)$ for a finitely generated field $k$; in that case $S$ is simply a smooth $k$-variety.
\item The case where $\Bb$ is an open subscheme of $\Spec(R)$ and $R$ is the ring of integers in a number field; in that case $S$ {{} is sometimes called \sl an arithmetic scheme.\rm} 
\end{enumerate}
 
 The following lemma is inspired by Drinfeld's {\sl ``conventional formulation of Hilbertianity''} in \cite[Section A.1.]{Drin}.  

\begin{lemm}  \label{connlemm}
If $K$ is Hilbertian, then for every dense open subscheme $U$ of $S$ and every finite \'etale 
morphism $p: X\to U$ there exists a connected subscheme $T$ of $U$ with the following properties.
\begin{enumerate}
\item[(a)] $T$ is a subscheme of $(f|U)^{-1}(\Bb')$ for some dense open subscheme $Z'$ of $\Bb$.
\item[(b)]  $f|T:T\to \Bb'$ is finite and \'etale.
\item[(c)] $p^{-1}(T)=X\times_S T$ is connected. 
\end{enumerate}
\end{lemm}

{\em Proof.} As $f$ is smooth, after replacing $U$ by a smaller dense open set and after replacing $X$ accordingly, there exists an \'etale 
$Z$-morphism
$g: U\to \Bb\times \Aa_d$ (cf. \cite[Expos\'e II, §1]{SGA1}). For any dense open subscheme $\Bb'$ of $\Bb$ we can consider the following commutative diagram of schemes
$$\xymatrix{
X_K\ar[r]^{p_K}\ar[d] & U_K\ar[r]^{g_K}\ar[d]& \Aa_{d, K}\ar[r]\ar[d] & Spec(K)\ar[d]\\
X'\ar[r]^{p'} & U'\ar[r]^{g'}& \Bb'\times \Aa_{d}\ar[r] & \Bb',
}$$
where $U'=U\times_\Bb \Bb'$, $X'=X\times_\Bb \Bb'$ and $g'$ and $p'$ are the restrictions of $g$ and $p$ respectively. 
As $K$ is Hilbertian there is a point $a\in\Aa_d(K)$ such that $(g_K\circ p_K)^{-1}(a)=\Spec(F)$  and $g_K^{-1}(a)=\Spec(E)$ where $E/K$ and $F/E$ 
are finite separable field extensions. For a suitable choice of $\Bb'$ the closed immersion $a: \Spec(K)\to \Aa_{d,K}$ extends to a closed subscheme
$Y$ of $\Bb'\times \Aa_d$ (cf. \cite[8.8.2 and 8.10.5]{EGAIV3}). Put $T:=g^{-1}(Y)$ and $X_T:=p^{-1}(T)$. After replacing $\Bb'$ by one of its dense open
subschemes we can assume that the maps  $p|X_T: X_T\to T$ and $g: T\to Y$ are finite and \'etale, because the corresponding maps on the generic fibres 
$$(g_K\circ p_K)^{-1}(a)=\Spec(F)\to g_K^{-1}(a)=\Spec(E)\to \Spec(K)$$
are finite and \'etale. The closed subscheme $T:=g^{-1}(Y)$ of $U'$ is connected because it is finite and \'etale over $\Bb'$ and its generic fibre 
$T\times_{\Bb'} \Spec(K)=g_K^{-1}(a)=\Spec(E)$ is connected. Likewise the closed subscheme $X_T$ of $X'$ is connected, because it is finite and \'etale over $\Bb'$ and its generic fibre $X_T\times_{\Bb'}\Spec(K)=\Spec(F)$ is connected.
\hfill $\Box$

\begin{coro}\label{rhogencoro}
Let $G$ be a profinite group and assume that the Frattini subgroup $\Phi(G)$ of $G$ is open in $G$. Let
$\rho: \pi_1(S)\to G$ be a group homomorphism. 
For every dense open subscheme $U$ of $S$ and every $m\in \{1,\cdots, d\}$ there exists a $T\in \Sm_m(S/\Bb)$ such that $T\subset U$ and such 
that $T$ is 
$\rho$-generic.
If $K$ is Hilbertian, then the same holds true for $m=0$. 
\end{coro}

{\em Proof.} We consider the homomorphism $\rho': \pi_1(S)\buildrel\rho\over\longrightarrow G\to G/\Phi(G).$ The image of $\rho'$ is finite because
$\Phi(G)$ is open in $G$. By the Frattini property, a subscheme $T\in  \Sm_m(S/\Bb)$ is $\rho$-generic if and only if it is $\rho'$-generic. In the proof of the corollary we can hence assume that $G$ is finite. 

{\bf Case A}:  Assume that $K$ is Hilbertian and $m=0$.  Let $X$ be the finite \'etale cover of $S$ corresponding to the kernel of $\rho$. 
By Lemma \ref{connlemm} there exists $T\in\Sm_m(S/\Bb)$ such that $T\subset U$ and such that $T\times_S X$  is connected. This
implies that $T$ is $\rho$-generic. 

{\bf Case B}: Assume that $m\in\{1,2,\cdots, d\}$ (and $K$ arbitrary). 
As $f: S\to \Bb$ is smooth we can, after replacing $U$ by a smaller open set, assume that there exists an \'etale $\Bb$-morphism $U\to \Bb\times \Aa_d$. Composing with an appropriate projection we get a smooth morphism
$$S\to \Bb\times \Aa_d\to \Bb\times \Aa_m.$$
Note that $\Sm_0(S/\Bb\times \Aa_m)\subset \Sm_m(S/\Bb)$. 
The function field $K(x_1,\cdots,x_m)$ of $\Bb\times \Aa_m$ is Hilbertian (even if $K$ is not). We can thus apply Case A with $\Bb$ replaced by
$\Bb\times\Aa_m$ to finish up the proof in Case B.\hfill $\Box$

\begin{rema} If $G$ is a compact subgroup of $\GL_n(\Qq_\ell)$, then $\Phi(G)$ is open in $G$ 
{{} (cf. 
\cite[Thm. 8.33]{Dix99}, \cite[§10.6]{Ser89}).} Hence, the above Corollary
\ref{rhogencoro} can be applied to $\ell$-adic representations of $\pi_1(S[\ell^{-1}])$, e.g., to $\rho_{A\times B,\ell}$. 
\end{rema}

 {\em From now on until the end of this section, let $\AAA$ and $\BBB$ be abelian schemes over $S$.  For $T\in\Sm_m(S/\Bb)$ we denote (as in the previous section) by $A_T/R(T)$ (resp. by $B_T/R(T)$) the generic fibre of 
 $\AAA_T\to T$ (resp. of $\BBB_T\to T$). Finally, for $\ell\in\Ll(S)$ 
 we define $\rho_\ell:=\rho_{A\times B, \ell}$.}\par\medskip
 
 The following lemma is at the core of the rest of the arguments of this section.

\begin{lemm}\label{fingenhomspezplus}
Let $m\in\{0,1,\cdots, d\}$ and $T\in\Sm_m(S/\Bb)$. 
Let $\Delta_1, \Delta_2$ be subsets of $\Nn$. 
If there exists  an $R(T)$-homomorphism $f: A_T\to B_T$ 
(resp. $\ol{R(T)}$-homomorphism $f: A_{T, \ol{R(T)}} \to B_{T,\ol{R(T)}}$) such that $\dim(\ker(f))\in \Delta_1$ and $\dim(\im(f))\in \Delta_2$ and if $T$ is 
$\rho_\ell$-generic, then
there exists an $R(S)$-homomorphism $F: A\to B$ (resp. $\ol{R(S)}$-homomorphism $F: A_{\ol{R(S)}}\to B_{\ol{R(S)}}$) such that $\dim(\ker(F))\in \Delta_1$ and $\dim(\im(F))\in \Delta_2$. 
\end{lemm}

{\em Proof.} Let  $f: A_T\to B_T$ be an $R(T)$-homomorphism such that we have $\dim(\ker(f)){\in}\Delta_1$ and $\dim(\im(f)){\in} \Delta_2,$ and assume 
that $T$ is $\rho_\ell$-generic. By Lemma \ref{fingenhomspez} the specialization maps
\begin{align*}
&r:=r_{\AAA,\BBB,T,S}: \Hom_{R(S)}(A,B) \to \Hom_{R(T)}(A_T, B_T),\\
&\ol{r}:=\ol{r}_{\AAA,\BBB,T,S}: \Hom_{\ol{R(S)}}(A_{\ol{R(S)}},B_{\ol{R(S)}}) \to \Hom_{\ol{R(T)}}(A_{T,\ol{R(T)}}, B_{T,\ol{R(T)}}),
\end{align*}
are injective with finite cokernels $C=\coker(r)$ and $\ol{C}=\coker(\ol{r})$. Let $s=|C|$. 
Then $s\circ f$ lies in the image of $r$. Thus there exists a $R(S)$-homomorphism $F: A\to B$ such that
$r(F)=s\circ f$. Using Lemma \ref{specializehoms} once more, we see that
\begin{align*}
&\dim(\ker(F))=\dim(\ker(s\circ f))=\dim(\ker(f))\in\Delta_1,\\
&\dim(\im(F))=\dim(\im(s\circ f))=\dim(\im(f))\in\Delta_2.
\end{align*}
The proof of the resp case can be carried out in a completely analoguous way using $\ol{r}$ instead of $r$ and taking Remark
\ref{rema:spec} into account. \hfill $\Box$

\medskip

We are ready for our first local-global statement. 

\begin{thm} \label{hilbertmain0}  Let $U$ be a dense open subscheme of $S$. 
Let $m\in \{0,1,\cdots, d\}$.  Assume that $K$ is Hilbertian or that $m\ge 1$. Let $\Delta_1$ and $\Delta_2$ be subsets of $\Nn$. 
The following are equivalent:
\begin{enumerate}
\item[(a)] There exists an $R(S)$-homomorphism $F: A\to B$ (resp. $\ol{R(S)}$-homomor\-phism $F: A_{\ol{R(S)}}\to B_{\ol{R(S)}}$)
such that $\dim(\ker(F))\in \Delta_1$ and one has $\dim(\im(F))\in\Delta_2$. 
\item[(b)] For every  $T\in \Sm_m(S/\Bb)$ with $T\subset U$ there exists a $R(T)$-homomor\-phism $f: A_T\to B_T$ 
(resp. $\ol{R(T)}$-homomor\-phism $f: A_{T, \ol{R(T)}} \to B_{T,\ol{R(T)}}$)
such that
$\dim(\ker(f))\in \Delta_1$ and one has $\dim(\im(f))\in\Delta_2$. 
\end{enumerate}
\end{thm}

{\em Proof.}
Assume (a) holds true and let $T\in \Sm_m(S/\Bb)$. There exists $\ell\in\Ll(T)$
by Lemma \ref{ellex}. 
The existence of $f$ as in (b) is then immediate from the mere existence of the specialization maps
$r_{\AAA,\BBB,T,S}$ and $\ol{r}_{\AAA,\BBB,T,S}$ (cf.  definitions of maps (\ref{spz1}) and (\ref{rbar}))  
plus the fact that they ``respect dimension data'' (cf. 
Lemma \ref{specializehoms} and Remark \ref{rema:spec}). This proves (a)$\Rightarrow$(b). 

We now prove the implication (b)$\Rightarrow$(a).
Assume (b) is true.  Choose $\ell\in\Ll(U)$ (cf. Lemma \ref{ellex}). By Corollary 
\ref{rhogencoro}  (applied with $S[\ell^{-1}]$ instead of $S$ and with $U[\ell^{-1}]$ instead of $U$))
there exists $T\in \Sm_m(S[\ell^{-1}]/\Bb)$ such that $T\subset U[\ell^{-1}]$ and such that $T$ is 
$\rho_\ell$-generic.
By (b) there exists a $R(T)$-homomorphism $f: A_T\to B_T$ 
(resp. $\ol{R(T)}$-homomorphism $f: A_{T, \ol{R(T)}} \to B_{T,\ol{R(T)}}$) such that
$\dim(\ker(f))\in \Delta_1$ and $\dim(\im(f))\in\Delta_2$. Now (a) follows by Lemma \ref{fingenhomspezplus}.\hfill $\Box$

\medskip

\noindent
{\em Proof of Theorem \ref{thmA}.} Apply Theorem \ref{hilbertmain0} with $\Delta_1{=}\{0\}$ and $\Delta_2{=}\{\dim(B)\}$ (resp. with $\Delta_1{=}\Nn$ and $\Delta_2{=}\{\dim(B)\}$, resp. with
$\Delta_1{=}\{0,1,{\cdots},\dim(A){-}1\}$ and $\Delta_2{=}\Nn$, resp. with $\Delta_1{=}\{\kappa\}$ and $\Delta_2{=}\Nn$).  \hfill $\Box$

\begin{coro} \label{hilbertmain2} Let $U$ be a dense open subscheme of $S$. 
Let $m\in \{0,1,{\cdots},d\}$.  Assume that $K$ is Hilbertian or that $m\ge 1$. 
\begin{enumerate}
\item[(a)] The following are equivalent:
\begin{enumerate}
\item[(i)] $A$ is not a simple $R(S)$-variety.
\item[(ii)] For every $T\in \Sm_m(S/\Bb)$ with $T\subset U$ 
the fibre $A_T$ is not a simple $R(T)$-variety.
\end{enumerate}
\item[(b)] The following are equivalent:
\begin{enumerate}
\item[(i)] $A_{\ol{R(S)}}$ is not a simple $\ol{R(S)}$-variety.
\item[(ii)] For every
 $T\in \Sm_m(S/\Bb)$ with $T\subset U$ 
the fibre $A_{T, \ol{R(T)}}$ is not a simple $\ol{R(T)}$-variety
\end{enumerate}
\end{enumerate}
\end{coro}

{\em Proof.} Note that $A$  is non-simple if and only if there exists $\kappa\in \{1,2,\cdots, $ 
$\dim(A){-}1\}$ and a homomorphism $A\to A$ with $\kappa$-dimensional kernel. Thus the corollary is a formal consequence of 
Theorem \ref{hilbertmain0}.\hfill $\Box$

\begin{thm} \label{hilbertmain3} Let $U$ be a dense open subscheme of $S$. 
Let $m{\in} \{0,1,\cdots,  d\}$.  Assume that $K$ is Hilbertian or that $m\ge 1$. 
The following are equivalent:
\begin{enumerate}
\item[(a)] $A$ is a quadratic isogeny twist of $B$
\item[(b)] For every $T\in \Sm_m(S/\Bb)$ with $T\subset U$ 
the abelian variety $A_T$ is a quadratic isogeny twist of $B_T$.
\end{enumerate}
The implication (a)$\Rightarrow$(b) holds true also in the case where $\Bb=\Spec(k)$ with a finite field $k$ and $m=0$. 
\end{thm}

{\em Proof. } We prove the implication (a)$\Rightarrow$(b): Assume that $A$ is a quadratic isogeny twist of $B$ and let $T\in \Sm_m(S/\Bb)$.  Assume that $T\subset U$.  Choose 
$\ell\in\Ll(T)$ (cf. Lemma \ref{ellex}).
It follows by Lemma \ref{TwistLemm},  that there exists a $\chi\in \Hom(\Gal(R(S)), \{\pm Id_B\})$ such that $\rho_{A, \ell}\cong \chi\otimes \rho_{B, \ell}$. {{}For every 
$\sigma\in \ker(\Gal(R(S))\to \pi_1(U[\ell^{-1}]))$ we have $\rho_{A, \ell}(\sigma)=Id_{T_\ell(A)}$ and $\rho_{B,\ell}(\sigma)=Id_{T_\ell(B)}$, hence $\chi(\sigma)=+Id_B$, i.e., $\chi$ factors through 
$\pi_1(U[\ell^{-1}])$.  
Now we have 
\begin{align*}
\rho_{A_T, \ell}\cong u^*\rho_{A, \ell}=u^*(\chi\otimes \rho_{B, \ell})=u^*\chi\otimes u^*\rho_{B, \ell}\cong \chi_T\otimes \rho_{B_T, \ell}
\end{align*}
where $u{:}T{\to} S$ is the embedding and $\chi_T{:} \pi_1(T[\ell^{-1}]){\to} \{\pm Id_{B_T}\}$ is a quadra\-tic character. By Lemma \ref{TwistLemm} $A_T$ is 
a quadratic isogeny twist of $B_T$.

We now prove the other implication (b)$\Rightarrow$(a): Assume that for every  $T\in \Sm_m(S/\Bb)$ with $T\subset U$ 
the fibre $A_T$ is a quadratic isogeny twist of $B_T$. Let $\ell\in\Ll(U)$  (cf. Lemma \ref{ellex}).
By Corollary \ref{rhogencoro} there exists $T\in \Sm_m(S/\Bb)$ with $T\subset U[\ell^{-1}]$ 
such that $T$ is 
$\rho_{A\times B, \ell}$-generic.  Let 
$$\KKK=\ker(\rho_{A\times B, \ell}: \pi_1(S[\ell^{-1}])\to GL_{T_\ell(A\times B)}(\Zz_\ell))$$ 
and $\KKK_T=\ker(\rho_{A\times B, \ell}\circ u_*)$. Then $u_*: \pi_1(T)\to \pi_1(S[\ell^{-1}])$ induces an isomorphism $\pi_1(T)/\KKK_T\to \pi_1(S[\ell^{-1}])/\KKK$.  By Lemma \ref{TwistLemm} there exists $\chi_T\in \Hom(\pi_1(T), \{\pm Id_{B_T}\})$ such that $\rho_{A_T, \ell}\cong 
\chi_T\otimes \rho_{B_T, \ell}$. It follows that $\chi_T(\sigma)=+Id_{B_T},$ for every $\sigma\in \KKK_T$, because $\sigma\in \ker(\rho_{A_T, \ell})\cap \ker(\rho_{B_T, \ell})$. Hence $\chi_T$ factors through $\pi_1(T)/\KKK_T$ and thus induces a quadratic cha\-racter $\chi$ of $\pi_1(S[\ell^{-1}])/\KKK$ and thus of 
$\pi_1(S[\ell^{-1}]),$ such that 
$\rho_{A, \ell}{\cong} \chi\otimes \rho_{B, \ell}$. Thus (a) holds true by Lemma \ref{TwistLemm}. \hfill $\Box$}



\section{Specialization to a finite field}
In Theorem \ref{hilbertmain0}, Corollary \ref{hilbertmain2} and Theorem \ref{hilbertmain3} the case when $k$ is finite and $m{=}0$ has been excluded. This case is more difficult
and cannot be handled with Hilbertianity alone. In this section we generalize parts of  {{}results \ref{hilbertmain0}, \ref{hilbertmain2}, \ref{hilbertmain3}}  (under additional assumptions) to this ``critical'' case by using recent results of Khare-Larsen \cite{KL20} and Fit\'e \cite{Fite21}. 

\begin{defi}\label{MWC} Let $A$ be an abelian variety over a finitely generated field $K$ and let $\ell\neq \chara(K)$ be a rational prime. Let $\uG_\ell$ be the connected component of the 
Zariski closure of $\rho_{A, \ell}(\Gal(K))$. 
We say that $A$ satisfies condition $MWC(A)$ if for all rational primes $\ell\neq \chara(K)$ the action of $\uG_{\ell, \ol{\Qq}_\ell}$ on each irreducible factor of the representation
$V_\ell(A)\otimes_{\Qq_\ell}\ol{\Qq}_\ell$ is minuscule in the sense of Bourbaki \cite{bourbaki-lie} cf. \cite[page 1]{KL20}. 
\end{defi}

\begin{rema} \label{MWCrema}
Let $A$ be an abelian variety over a finitely generated field $K$ satisfying $MWC(A)$. Let $B$
be an abelian variety. 
\begin{enumerate}
\item[(a)] For every finite extension $K'/K$ the abelian variety $A_{K'}$ satisfies $MWC(A_{K'})$ because
the connected component of the 
Zariski closure of $\rho_{A, \ell}(\Gal(K'))$ agrees with  the connected component of the 
Zariski closure of $\rho_{A, \ell}(\Gal(K))$ as $\Gal(K')$ is of finite index in $\Gal(K)$. 
\item[(b)] If $f: A\to B$ is a surjective $K$-homomorphism, then $B$ satisfies $MWC(B)$, because $V_\ell(f): V_\ell(A)\to
V_\ell(B)$ is  a surjective homomorphism of representations. 
\item[(c)] If $B$ is an abelian subvariety of $A$, then there exists surjective homomorphism
$A\to B$ by \cite[Prop. 12.1]{milneav} and thus $B$ satisfies $MWC(B)$ by (b).
\end{enumerate}
\end{rema}

\begin{rema}\label{RMWC}
If $K$ is a number field, then $MWC(A)$ holds true for every abelian variety $A/K$ due to a result of Pink \cite[Cor. 5.11]{Pink98} which proves Zarhin's {{}\sl minuscule weights conjecture\rm} \cite[Conjecture 0.4]{Zar85} in the number field case.   The same holds true if $K$ is a finitely generated field of
characteristic zero; using Hilbertianity one can easily reduce to the number field case.
According to \cite[Conjecture 0.4]{Zar85} every abelian variety over a global field of any characteristic should satisfy $MWC(A)$. Zarhin \cite[4.2.1]{Zar85} proved his conjecture for $K$ a global field of positive characteristic and $A/K$ an ordinary abelian variety. In a recent preprint \cite{CT20} Cadoret and Tamagawa formulated an anologue of Zarhin's conjecture over finitely generated fields and checked that $MWC(A)$ holds true for another case of $A$ over $K$ finitely generated of positive characteristic, see \cite[Cor. 6.3.2.3]{CT20}. On the other hand, the analogue of the conjecture of Zarhin in positive characteristic seems to fail in general, as indicated by a result of B\" ultel (cf. \cite[Thm. 1.2]{Bue},
\cite[Comments and examples, p. 637]{Bue}) who constructed an abelian variety of dimension $56$ over a finitely generated field of positive characteristic with a simple factor of exceptional type $G_2$ (there are no minuscule weights in this case !) in its $\ell$-adic monodromy 
group. \footnote{For more examples of abelian varieties of this type the reader can consult the pre\-print \cite{Bue2} which appeared in the arXives after our paper had been submitted for publication.}
\end{rema}

\noindent
The following is a global function field analogue of a result of  Khare and Larsen \cite[Thm. 1]{KL20}. We include the proof for the reader's convenience. 

\begin{prop} \label{KL-lemm} Let $S$ be a smooth curve over a finite
field $k$. Let $\AAA, \BBB$ be abelian schemes over $S$ with
generic fibers $A$ and $B$ respectively. Assume that $A$ satisfies 
$MWC(A)$ and $B$ satisfies $MWC(B)$. If the set of all closed points $s\in S$ such that there exists a  non-zero   $\ol{k(s)}$-homomorphism 
$A_{s, \ol{k(s)}}\to B_{s, \ol{k(s)}}$ has Dirichlet density $1$, then there exists a non-zero $\ol{k(S)}$-homomorphism $A_{\ol{k(S)}}\to B_{\ol{k(S)}}$.
\end{prop}

{\em Proof.} Let $K=R(S)$. Fix one rational prime $\ell\neq \chara(K)$. 
After replacing $S$ by one of its connected finite \'etale covers, we can assume that the Zariski closures $\uG_A$ of $\rho_{A, \ell}(\Gal(K))$, $\uG_B$ of $\rho_{B, \ell}(\Gal(K))$ and $\uG$ of $\rho_{A\times B, \ell}(\Gal(K))$ are connected. The action of $\uG_{\ol{\Qq}_\ell}$ on $V_A:=V_\ell(A)\otimes_{\Qq_\ell} \ol{\Qq}_\ell$ and on $V_B:=V_\ell(B)\otimes_{\Qq_\ell} \ol{\Qq}_\ell$ is minuscule. By \cite[Thm. 1.2]{LP97}
 and \cite{Ser81} the set of all closed points $s\in S$ such that the Frobenius element $Fr_s\in \uG(\Qq_\ell)$ generates a Zariski dense subgroup $\langle Fr_s\rangle$
 of a maximal torus of $\uG$ has positive density. Thus we can choose a closed point $s\in S$ 
 such that $\langle Fr_s\rangle$ is a Zariski dense subgroup  of a maximal torus $\uT$ of $\uG$ and such that
 there exists a non-zero homomorphism $f: A_{s, \ol{k(s)}}\to B_{s, \ol{k(s)}}$. Furthermore $f$ is defined over a finite extension $k'$ of 
 $k(s)$. Now $|k'|=|k(s)|^m,$ for a natural number $m$ and thus $f$ is fixed under the action of $Fr_s^m$, and $Fr_s^m$ also generates a Zariski dense subgroup of $\uT$. It follows that 
 the map $V_\ell(f)\otimes\ol{\Qq}_\ell: V_A\to V_B$ induced by $f$ is fixed by 
 the torus $\uT$.  Moreover  $V_\ell(f)\otimes\ol{\Qq}_\ell$ is a non-zero element of $\Hom_\uT(V_A, V_B)$. Thus
 $$\dim(\Hom_\uT(V_A, V_B))>0.$$
 Now $\uG$ is a reductive group because $\rho_{A\times B, \ell}$ is semisimple. 
 By our assumption $A$ satisfies $MWC(A)$ and $B$ satisfies $MWC(B)$. Thus we can apply 
 \cite[Prop. 2]{KL20} to conclude that 
 $$\dim(\Hom_\uG(V_A, V_B))>0.$$
 From $\Hom_\uG(V_A, V_B)\cong \Hom_K(A,B)\otimes {\ol{\Qq}_\ell}$ and the fact that $\Hom_K(A,B)$ is $\Zz$-free  it follows that $\Hom_K(A,B)$ is non-zero, as desired.\hfill $\Box$ 
\medskip

We can now treat higher dimensional $S$ {{}by} combining Proposition \ref{KL-lemm} with the results of the previous section.

\begin{thm} \label{KL-thm} Let $S$ be a smooth variety over a finite
field $k$. Let $\AAA, \BBB$ be abelian schemes over $S$ with
generic fibers $A$ and $B$ respectively. 
Assume that $A$ satisfies 
$MWC(A)$ and $B$ satisfies $MWC(B)$. 
The following are equivalent:
\begin{enumerate}
\item[(a)] There exists a non-zero $\ol{k(S)}$-homomorphism $A_{\ol{k(S)}}\to B_{\ol{k(S)}}$.
\item[(b)] For every closed point $s\in S$ there exists a  non-zero  $\ol{k(s)}$-homomor\-phism 
$A_{s, \ol{k(s)}}\to B_{s, \ol{k(s)}}$. 
\end{enumerate}
\end{thm}

{\em Proof.} After replacing $S$ by one of its connected finite \'etale covers we can assume that 
$$\Hom_{k(S)}(A, B)=\Hom_{\ol{k(S)}}(A_{\ol{k(S)}}, B_{\ol{k(S)}}).$$
The implication (a)$\Rightarrow$(b) is immediate from Lemma
\ref{specializehoms} and Remark \ref{rema:spec}.
We prove the other 
implication (b)$\Rightarrow$(a):  Let $K=R(S)$ and let $\ell\neq\chara(K)$ be a rational prime.
By Corollary \ref{rhogencoro} there exists a $T\in \Sm_1(S/k)$ of $S$ such that $T$ is $\rho_{A\times B, \ell}$-generic. Then $T$ is automatically
$\rho_{A, \ell}$-generic and $\rho_{B, \ell}$-generic. In particular, $\rho_{A, \ell}(\Gal(K))^{\Zar}=\rho_{A_T, \ell}(\Gal(k(T)))^{\Zar}$. 
In particular, $A_T/k(T)$ satisfies $MWC(A_T)$. Similarly $B_T/k(T)$ satisfies $MWC(B_T)$. By (b), for 
every closed point $t\in T$ there exists a non-zero homomorphism $A_{\ol{k(t)}}\to B_{\ol{k(t)}}$. By 
Proposition \ref{KL-lemm} there exists a non-zero homomorphism $A_{T, \ol{k(T)}}\to B_{T, \ol{k(T)}}$. 
Now (a) follows by Lemma \ref{fingenhomspezplus}.\hfill $\Box$ 
\medskip

We can upgrade Theorem \ref{KL-thm} a bit, to treat not only non-zero homomorphisms but also other important classes of homomorphisms, much in the spirit of Corollary \ref{thmB}.

\begin{thm} \label{KL-thm2} Let $S$ be a smooth variety over a finite
field $k$. Let $\AAA, \BBB$ be abelian schemes over $S$ with
generic fibers $A$ and $B$ respectively. 
Assume that $A$ satisfies 
$MWC(A)$ and $B$ satisfies $MWC(B)$. 
The following are equivalent:
\begin{enumerate}
\item[(a)] There exists a surjective $\ol{k(S)}$-homomorphism (resp. $\ol{k(S)}$-isogeny) $A_{\ol{k(S)}}\to B_{\ol{k(S)}}$.
\item[(b)] For every closed point $s\in S$ there exists a  surjective  $\ol{k(s)}$-homomor\-phism (resp. $\ol{k(s)}$-isogeny)
$A_{s, \ol{k(s)}}\to B_{s, \ol{k(s)}}$. 
\end{enumerate}
\end{thm}

{\em Proof.} The implication (a)$\Rightarrow$(b) is immediate from Lemma \ref{specializehoms}.

We prove the implication (b)$\Rightarrow$ (a) in case ``surjective''. 
After replacing $S$ by one of its connected finite \'etale covers we can assume that 
\begin{align}\label{eq:End}
&\Hom_{k(S)}(A, B)=\Hom_{\ol{k(S)}}(A_{\ol{k(S)}}, B_{\ol{k(S)}}),\\\notag
&\End_{k(S)}(A)=\End_{\ol{k(S)}}(A_{\ol{k(S)}}),\\\notag
&\End_{k(S)}(B)=\End_{\ol{k(S)}}(B_{\ol{k(S)}}).
\end{align}
By the Poincar\'e reducibility theorem (cf. \cite[Prop. 12.1]{milneav} and the passage below) $A$ (resp. $B$) is $k(S)$-isogenous to $\prod_{i\in I} A_i^{n_i}$
(resp. $\prod_{j\in J} B_j^{m_i}$) where $I$ and $J$ are finite sets, 
the $A_i$ ($i\in I$) are mutually not $k(S)$-isogenous 
simple abelian varieties over $k(S)$ and the $B_j$ ($j\in J$)) are mutually not $k(S)$-isogenous 
simple abelian varieties over $k(S)$. 
We consider the commutative diagram 
$$\xymatrix{
\End^0_{\ol{k(S)}}(A_{\ol{k(S)}})\ar[r] & \prod_{i,j\in I} 
\Hom^0_{\ol{k(S)}}(A_{i, \ol{k(S)}}, A_{j, \ol{k(S)}})^{n_i\times n_j}\\
\End^0_{k(S)}(A_{k(S)})\ar[r] \ar[u]& \prod_{i,j\in I} 
\Hom^0_{k(S)}(A_{i}, A_{j})^{n_i\times n_j}\ar[u]_{\prod_{i,j\in I} s_{ij}^{n_i\times n_j}}\\
}
$$
where the $$s_{ij}: \Hom^0_{k(S)}(A_i, A_j)\to \Hom^0_{\ol{k(S)}}(A_{i, \ol{k(S)}}, A_{j, \ol{k(S)}})$$
are the canonical maps. 
The horizontal maps in the diagram are bijective. The left hand vertical map is bijective by 
\eqref{eq:End}. Thus the $s_{ij}$ are bijective, too. 
 For $i\neq j$ we have $A_i\not\simeq A_j$, hence 
 $$\Hom^0_{\ol{k(S)}}(A_{i, \ol{k(S)}}, A_{j, \ol{k(S)}})=\Hom^0_{k(S)}(A_{i}, A_{j})=0,$$
and thus the 
 $A_{i, \ol{k(S)}}$ are mutually non-isogenous over $\ol{k(S)}$. As $A_i$ is simple it follows that
  $\End^0_{\ol{k(S)}}(A_{i, \ol{k(S)}})=\End^0_{k(S)}(A_{i})$ is a division algebra over $\Qq$ and thus the
  $A_{i, \ol{k(S)}}$ are simple. 
Likewise the 
$B_{j, \ol{k(S)}}$ are 
simple and mutually non-isogenous over $\ol{k(S)}$. Note that $MWC(A_i)$ and $MWC(B_j)$ holds true (cf. Remark \ref{MWCrema}).

After replacing $S$ by one of its dense open subschemes each $A_i$ (resp. $B_j$) extends to an abelian scheme $\AAA_i$ (resp. $\BBB_j$) over $S$. It is then clear that for every closed point 
$s\in S$ the fibre $A_{s, \ol{k(s)}}$ is $\ol{k(s)}$-isogenous to $\prod_{i\in I} A_{i, s, \ol{k(s)}}^{n_i}$ and
 $B_{s, \ol{k(s)}}$ is $\ol{k(s)}$-isogenous to $\prod_{j\in J} B_{j, s, \ol{k(s)}}^{m_j}$. In what follows we can thus asssume, without loss of generality,  that $A=\prod_{i\in I} A_i^{n_i}$ and that $B=\prod_{j\in J} B_j^{m_j}.$ From
(b) we know that for every closed point $s\in S$ there exists a  surjective  $\ol{k(s)}$-homomorphism 
$A_{s, \ol{k(s)}}\to B_{s, \ol{k(s)}}$ and thus for every $j\in J$ a surjective homomorphism
$A_{s, \ol{k(s)}}\to B_{j, s, \ol{k(s)}}$.  By Theorem \ref{KL-thm} there exists for every $j\in J$ a non-zero homomorphism 
$A_{\ol{k(S)}}\to B_{j, \ol{k(S)}}$. 
Thus, for every $j\in J$, there must be a unique $\nu(j)\in I$, such that there exists a non-zero $\ol{k(S)}$-homomorphism
$h_j: A_{\nu(j), \ol{k(S)}}\to B_{j,\ol{k(S)}}$ and $h_j$ must be a $\ol{k(S)}$-isogeny, because 
$A_{\nu(j), \ol{k(S)}}$ and $B_{j,\ol{k(S)}}$ are both simple. 

We now show that $n_{\nu(j)}\ge m_j,$ for all $j\in J$.
Let $j\in J$ be arbitrary. Let $A'_j=\prod_{i\in I\setminus \{\nu(j)\}} A_i^{n_i}$. Then for every $j\in J$ we have $A=A_{\nu(j)}^{n_{\nu(j)}}\times A_j'$. There is no
non-zero 
$\ol{k(S)}$-homomorphism $A_{j, \ol{k(S)}}'\to B_{j,\ol{k(S)}}$, because the $A_{i, \ol{k(S)}}$ are mutually 
non-isogenous over $\ol{k(S)}$ and $A_{\nu(j), \ol{k(S)}}\simeq B_{j, \ol{k(S)}}$. 
The contraposition of Theorem \ref{KL-thm} implies that 
there exists a closed point $s\in S$ such that there is no non-zero $\ol{k(s)}$-homomorphism 
$A_{j, s, \ol{k(s)}}'\to B_{j, s, \ol{k(s)}} $. Furthermore, by (b), there is a surjective 
$\ol{k(s)}$-homomorphism  
$$A_{s, \ol{k(s)}}=A_{\nu(j), s, \ol{k(s)}}^{n_{\nu(j)}}\times A_{j, s, \ol{k(s)}}'\to B_{j, s, \ol{k(s)}}^{m_j}$$
which restricts to a surjective $\ol{k(s)}$-homomorphism  
$$A_{\nu(j), s, \ol{k(s)}}^{n_{\nu(j)}}\times \{0\}\to B_{j, s, \ol{k(s)}}^{m_j}.$$
As $\dim(A_{\nu(j), s, \ol{k(s)}})=\dim(B_{j, s, \ol{k(s)}})$, 
this finally shows that $n_{\nu(j)}\ge m_j$ as claimed above. 

This together with the existence of the isogeny $h_j$ implies the existence of a surjective $\ol{k(S)}$-homomorphism $f_j: A_{\nu(j), \ol{k(S)}}^{n_{\nu(j)}} \to B_{j, \ol{k(S)}}^{m_j},$ and these in turn induce 
a surjective $\ol{k(S)}$-homomorphism
\begin{align}\label{eq:prodhom}
\prod_{j\in J} f_j : \prod_{j\in J} A_{\nu(j), \ol{k(S)}}^{n_{\nu(j)}} \to \prod_{j\in J} B_{j, \ol{k(S)}}^{m_j}=
B_{\ol{k(S)}}.
\end{align}
Composing the surjective $\ol{k(S)}$-homomorphism $\eqref{eq:prodhom}$ with the 
projection $A_{\ol{k(S)}}\to  \prod_{j\in J} A_{\nu(j), \ol{k(S)}}^{n_{\nu(j)}}$ we obtain a 
$\ol{k(S)}$-epimorphism $A_{\ol{k(S)}}\to B_{\ol{k(S)}}$ as desired.

We next prove the implication (b)$\Rightarrow$(a) in case ``isogeny''. In that case we know by (b) that
for every closed point $s\in S$ there exists a  $\ol{k(s)}$-isogeny
$g_s: A_{s, \ol{k(s)}}\to B_{s, \ol{k(s)}}$. The above arguments imply the 
existence of a surjective $\ol{k(S)}$-homo\-morphism $g: A_{\ol{k(S)}}\to B_{\ol{k(S)}}$.  Let $s\in S$ be 
a closed point. Then, by the existence of the isogeny $g_s$, $$\dim(A)=\dim(A_{s, \ol{k(s)}})=\dim(B_{s, \ol{k(s)}})=\dim(B),$$ and thus
$g$ must automatically be an isogeny. \hfill $\Box$
\par\medskip

We shall now establish a global function field analogue of \cite[Cor. 2.7]{Fite21}, following proof {{}in \sl loc. cit.\rm} quite closely. For this we need the following theorem 
of Rajan \cite{Ra98}. Let $K$ be a global field and denote by $\Sigma_K$ the set of all non-archimedian discrete valuations of $K$. Let $\ell\neq \chara(K)$ be a rational prime and let $\rho_1: \Gal(K)\to \GL_n(\Qq_\ell)$ and $\rho_2: \Gal(K)\to \GL_n(\Qq_\ell)$ be continuous semisimple representations which are unramified outside a finite set $S\subset \Sigma_K$. Define
$$SM(\rho_1, \rho_2):=\{v\in \Sigma_K: Tr(\rho_1(Fr_v))=Tr(\rho_2(Fr_v))\}.$$

\begin{thm} (Rajan, \cite[Thm. 2]{Ra98}) \label{thm:rajan} Suppose that the Zariski closure $H_1$ of $\rho_1(\Gal(K))$ is connected and that the upper Dirichlet density of $SM(\rho_1, \rho_2)$ is 
 positive. Then the following holds true.
\begin{enumerate}
\item[(a)] There is a finite Galois extension $L/K$ such that $\rho_1|_{\Gal(L)}{\cong} \rho_2|_{\Gal(L)}$ and the connected component of the Zariski closure of $\rho_2(\Gal(K))$ is
conjugate to $H_1$.
\item[(b)] If $\rho_1$ is absolutely irreducible, then there is a Dirichlet character $\chi: \Gal(L/K)\to \Qq_\ell^\times$ of finite order such that $\rho_2\cong \chi\otimes_{\Qq_\ell} \rho_1$. 
\end{enumerate}
\end{thm}

{{}We have the following global function field analogue of \cite[Cor. 2.7]{Fite21}.} 

 \begin{prop} \label{Fite-lemm} Let $S$ be a smooth curve over a finite
field $k$. Let $\AAA, \BBB$ be abelian schemes over $S$ with
generic fibers $A$ and $B,$ respectively. Let $\ell\neq\chara(k)$ be a rational prime {{} and let $K=k(S)$.} 
Let $U$ be a dense open subscheme of $S$. 
Assume that $\End_{\ol{k(S)}}(A)=\End_{\ol{k(S)}}(B)=\Zz$ and that {{}the Zariski closures of $\rho_{A, \ell}(\Gal(K))$ and of 
$\rho_{B, \ell}(\Gal(K))$ are connected.}
If the set of all  closed points $u\in U$ such that 
$A_u$ is a quadratic isogeny twist of $B_u$ has Dirichlet density $1$, then $A$ is a quadratic isogeny twist of $B$. 
\end{prop}

{\em Proof.} Let $\Gamma$ be the set of all $u\in U$ such that 
$A_u$ is a quadratic isogeny twist of $B_u$. Assume that $\Gamma$ has Dirichlet density $1$. We claim that $SM(\rho_{\ell, A}, \rho_{\ell, B})$ has positive upper 
Dirichlet density. Consider the two sets $$\Gamma_{\pm}:=\{u\in\Gamma: Tr(\rho_{A, \ell}(Fr_u))=\pm Tr(\rho_{B, \ell}(Fr_u))\}.$$
For every $u \in\Gamma$ we have $Tr(\rho_{A, \ell}(Fr_u))\in \{\pm Tr(\rho_{B, \ell}(Fr_u))\}$ by Lemma \ref{TwistLemm} (implication (a)$\Rightarrow$(c) applied over $k(u)$). Thus $\Gamma=\Gamma_+\cup \Gamma_-$. Moreover $\Gamma_+\subset  SM(\rho_{\ell, A}, \rho_{\ell, B})$. {{}If the inclusion were wrong,} then it would follow that $\Gamma_+$ has upper Dirichlet density zero and that $\Gamma_-$ has Dirichlet density $1$. But then the Chebotarev density theorem would imply that $Tr(\rho_{A, \ell}(g))=-Tr(\rho_{B, \ell}(g))$ for all $g\in\Gal(K)$, which is obviously false for $g=Id$. Thus {{}the inclusion $\Gamma_+\subset  SM(\rho_{\ell, A}, \rho_{\ell, B})$} holds true. 

Our assumption $\End_{\ol{k(S)}}(A)=\End_{\ol{k(S)}}(B)=\Zz$ implies that $\rho_{A, \ell}$ and $\rho_{B, \ell}$ are absolutely irreducible. 

Now, because of the claim, Theorem \ref{thm:rajan} implies that there is a finite Galois extension $L/K$ and a character $\chi: \Gal(L/K)\to \Qq_\ell^\times$ of finite order such that $\rho_{A,\ell}=\chi\otimes \rho_{B, \ell}$. But now, by Theorem \ref{zf}, we have an isomorphism of $\Gal(L/K)$-modules
\begin{align*}
\Hom_L(A_L, B_L)\otimes\Qq_\ell&\cong V_\ell(A)^*\otimes_{\Qq_\ell} V_\ell(B)\cong\\
&\cong V_\ell(A)^*\otimes_{\Qq_\ell} V_\ell(A)\otimes \Qq_\ell(\chi)\cong\\
&\cong \End_L(A_L)\otimes \Qq_\ell(\chi)\cong \Qq_\ell(\chi).
\end{align*}
Thus $\Gal(L/K)$ acts on $\Hom_L(A_L, B_L)\cong \Zz$ via $\chi$ and this implies that $\chi$ is a {\em quadratic} character. Lemma \ref{TwistLemm} (implication
(b)$\Rightarrow$ (a) applied over $K$) now implies that $B$ is a quadratic isogeny twist of $A$. \hfill $\Box$
\medskip

Again one can eliminate the assumption $\dim(S)=1$ {{}by the Hilbertianity approach of the previous section.} 

 \begin{coro} \label{Fite-coro} Let $S$ be a smooth variety over a finite
field $k$. Let $\AAA, \BBB$ be abelian schemes over $S$ with
generic fibers $A$ and $B$ respectively. Let $U$ be a dense open subscheme of $S$. 
Assume that $End_{\ol{k(S)}}(A)=End_{\ol{k(S)}}(B)=\Zz$. 
The following are equivalent:
\begin{enumerate}
\item[(a)]  $A$ is a quadratic isogeny twist of $B$.
\item[(b)] For every closed point $u$ of $U$ the abelian variety
$A_u$ is a quadratic isogeny twist of $B_u$. 
\end{enumerate}
\end{coro}

{\em Proof.} The implication (a)$\Rightarrow$(b) is known from Theorem \ref{hilbertmain3}. We prove the other implication.
Assume that for every closed point $u$ of $U$ the abelian variety
$A_u$ is a quadratic isogeny twist of $B_u$. Then, for every smooth connected curve $T$ on $U$ the abelian variety $A_T$ is a quadratic isogeny twist of $B_T$ by Proposition \ref{Fite-lemm}. Thus Theorem \ref{hilbertmain3} implies (a). \hfill $\Box$


{{}\begin{rema}
It is easy to show that if $A$ and $B$ in Corollary \ref{Fite-coro} are elliptic curves with nontrivial endomorphisms, then the claim still holds. It is so, since then j-invariants $j(A)$ and $j(B)$ belong to $\ol{k},$ hence the curves are isotrivial cf. \cite[V.3.1 and Exerc.V.5.8]{Silv09}, so a twist between fibres at $u\in U$ extends to a twist between the curves. It is an interesting question, if Corollary \ref{Fite-coro} holds true for nonsimple abelian varieties, e.g.,  for products of mutually nonisogenous elliptic curves over $k(S).$  
\end{rema}}


\newpage

{\sc Wojciech Gajda,\\
Faculty of Mathematics and Computer Sciences\\
Adam Mickiewicz University,\\
Uniwersytetu Pozna\'nskiego 4\\ 
61-614 Pozna\'{n}, Poland}\\
\it E-mail: \tt gajda@amu.edu.pl
\bigskip\bigskip

{\sc Sebastian Petersen\\ 
Universit\"at Kassel\\
Fachbereich 10,\\
Wilhelmsh\"oher Allee 73\\
34-119 Kassel, Germany}\\
\it E-mail: \tt petersen@mathematik.uni-kassel.de
\end{document}